\documentclass[12pt]{article}
\usepackage{amsmath,amsfonts,amssymb,amsthm}
\usepackage{enumitem}
\usepackage{color}
\usepackage{authblk}

\newtheorem{proposition}{Proposition}
\newtheorem{theorem}{Theorem}

\theoremstyle{remark}
\newtheorem{remark}{Remark}
\newtheorem{example}{Example}

\newtheoremstyle{assumptionstyle}
  {\topsep}      
  {\topsep}      
  {\normalfont\leftskip=3em} 
  {-2.5em}         
  {\bfseries}    
  {.}            
  {1em}          
  {\thmnumber{#2}} 

\theoremstyle{assumptionstyle}
\newtheorem{assumption}{Assumption}

\DeclareMathOperator*{\olim}{\overline{\lim}}
\DeclareMathOperator*{\ulim}{\underline{\lim}}

\newcommand{\eqd}{\overset{\mathcal D}{=}}
\newcommand{\PP}{{\mathbf P}}
\newcommand{\E}{{\mathbf E}}

\newcommand{\nc}{\newcommand}

\nc{\I}{{\bf 1}}
\nc{\R}{{\mathbb R}}
\nc{\N}{{\mathbb N}}
\nc{\Z}{{\mathbb Z}}

\nc{\BP}{\mathbb{P}}
\nc{\BE}{\mathbb{E}}
\nc{\BQ}{\mathbb{Q}}
\nc{\BX}{\mathbb{X}}
\numberwithin{equation}{section}

\usepackage[style=authoryear,maxcitenames=2,]{biblatex}

\begin{document}
\renewcommand{\thefootnote}{\fnsymbol{footnote}}

\author[1,2]{Sergey Foss}
\author[3]{Peter W. Glynn}
\affil[1]{Heriot-Watt University}
\affil[2]{Sobolev Institute of Mathematics}
\affil[3]{Stanford University}
\title{
On Recurrence of the Infinite Server Queue}
\date{\today}
\maketitle

\begin{abstract}
This paper concerns the recurrence structure of the infinite server queue, as viewed through the prism of the maximum dater sequence, namely the time to drain the current work in the system as seen at arrival epochs. Despite the importance of this model in queueing theory, we are aware of no complete analysis of the stability behavior of this model, especially in settings in which either or both the inter arrival and service time distributions have infinite mean. In this paper, we fully develop the analog of the Loynes construction of the stationary version in the context of stationary ergodic inputs, extending earlier work of E.~Altman (2005), and then classify the Markov chain when the inputs are independent and identically distributed. This allows us to classify the chain, according to transience, recurrence in the sense of Harris, and positive recurrence in the sense of Harris. We further go on to develop tail asymptotics for the stationary distribution of the maximum dater sequence, when the service times have tails that are asymptotically exponential or Pareto, and we contrast the stability theory for the infinite server queue relative to that for the single server queue.
\end{abstract}

Keywords: Infinite server queue, Loynes lemma, infinite mean distributions

\section{Introduction}\label{S:In}

    In the analysis of simple single station queues, there are two models that take on special mathematical importance as fundamental models, namely the single-server queue with first-in first-out (FIFO) service  discipline and the infinite server queue. In the well-known Kendall notation for queues, the single server model is denoted as the $G/G/1$ queue, whereas the infinite server queue is denoted as the $G/G/\infty$ queue. The infinite server queue is of interest in its own right, and also adds useful insight into the behavior of many-server queues in which the number of servers is large. While a great deal is known about the recurrence behavior of the single-server queue, less is known about the recurrence of the infinite server queue. This paper therefore is focused on discussion of the recurrence structure for the infinite-server model.

   Since we wish to eventually draw comparisons between the recurrence behavior of the single-server queue and the infinite-server queue, we
    start by consideration of a fundamental process that is central to both models, namely the {\it virtual workload process} $W=(W(t): t\ge 0)$. The random variable (rv) $W(t)$ represents the amount of ``wall-clock" time the system would require to ``drain" all its work if no additional work were to arrive to the system after time $t$. If $T_n$ is the arrival time of the $n$'th customer, either of the two discrete-time sequences $(W(T_n):n \ge 0)$ and $(W(T_{n}{\scriptstyle -}): n\ge 0)$ can be used to study these models, where $W(t{\scriptstyle -})$ represents the left-limit of $W$ at $t$. For the single-server queue, the {\it waiting time} $W_n := W(T_{n}{\scriptstyle -})$ for the $n$'th customer is more typically studied, and satisfies the Lindley recursion
   $$
 W_{n+1}= \max (W_n+s_n - t_{n+1}, 0)
 $$
   for $n \ge 0$, where $s_n$ is the service time for customer $n$, and $t_{n+1} = T_{n+1} - T_n$ is the $(n+1)$'st inter-arrival time. On the other hand, for the infinite server queue, we will focus on the {\it maximum dater} sequence $X = (X_n: n \ge 0)$, where $X_n = W(T_n)$ for $n \ge 0$; see F.~Baccelli and S.~Foss (1995) for discussion of maximum daters. For the $G/G/\infty$ queue, this sequence satisfies the recursion
 \begin{equation}\label{SRS1}
 X_{n+1}= \max (X_n-t_{n+1},s_{n+1}),
 \end{equation}
  for $n \ge 0$, where $X_0=x$ is the initial workload. This recursion is also studied as a special case of the theory developed by F.~Comets, R.~Fernansez and P.~Ferrari (2002).
  
      Our main contributions in this paper relate to a reasonably complete development of the recurrence theory for the infinite-server queue, with the goal of contrasting this theory with the parallel development for the single-server queue. In Section 2, we study the theory in the setting of stationary ergodic input sequences. Our goal is to develop necessary and sufficient conditions under which the maximum dater for the infinite-server queue admits a stationary version; see Theorem 1 and Proposition 2. These results complement those of E.~Altman (2005), in which the stationary version is studied under the usual finite mean conditions on the inputs. Section 3 studies the recurrence problem in the setting of the $GI/GI/\infty$ queue. We obtain necessary and sufficient conditions in Theorem 2 for transience, Harris recurrence, and positive Harris recurrence, and then apply our criteria to a couple of infinite mean examples. Not surprisingly, the recurrence and transience depend upon how heavy is the service time tail, with null recurrence possible when (for example) the service time Pareto tail index is just on the cusp of having an infinite mean. In the case where one has Markov dependence and the service times have finite mean, the result is known; see e.g. Section 5 of C.M.~Goldie (1991).  S.~Popov (2025) develops related theory for the number-in-system process for $GI/GI/\infty$ when the arrivals are Poisson and the system is started empty. We conclude Section 3 by contrasting our necessary and sufficient condition for positive recurrence of $G/GI/\infty$ with the corresponding condition for the $GI/GI/1$ queue.

      Our final section develops exact asymptotics for the stationary version of the infinite-server queue's maximum dater sequence; see Theorem 3. Again, the tail asymptotics are contrasted with the corresponding asymptotics for $GI/GI/1$.

\section{The Loynes Construction for the Infinite Server Queue}\label{Loynes}

   In this section, we focus our discussion on the case in which the $(t_j, s_j)$'s are stationary and ergodic. This complements earlier work of E.~Altman (2005) in which the $s_j$'s and $t_j$'s  are assumed to have finite mean. More precisely, we assume the existence of a two-sided sequence:
   \begin{assumption} \label{a1}
      $((s_j, t_j) : -\infty < j < \infty)$ is a stationary ergodic sequence of positive random variables (rv's).
   \end{assumption}
   \noindent For $x \ge 0$, put $X_0(x) = x$ and
   \[ X_{n+1}(x) = \max(X_n(x) - t_{n+1}, s_{n+1}) \]
   for $n \ge 0$, so that $(X_n(x): n \ge 0)$ is the maximum dater sequence initialized at $x$. Let $T_0 = 0$ and $T_n = t_1 + \cdots + t_n$ be the arrival instant for customer $n$. We say that $Y_n$ converges in \textit{total variation} to $Y_\infty$ (and write $Y_n \overset{\text{tv}}{\to} Y_\infty$) if
   \[ \sup_A | \PP(Y_n \in A) - \PP(Y_\infty \in A) | \to 0 \]
   as $n \to \infty$. For $n \in \mathbb{Z}$, put $\tilde s_n = s_{-n}$ and $\tilde t_n = t_{-n}$, and for $0 \le n \le \infty$, set
   \[ \widetilde X_n = \max (0, \max_{1 \le j \le n} [\tilde s_j - \sum_{k=0}^{j-1} \tilde t_k]). \]
   For rv's $Y_1$ and $Y_2$, we write $Y_1 \eqd Y_2$ when $Y_1$ and $Y_2$ have the same distribution.

   \begin{proposition}
      ~\label{p1}
      \begin{enumerate}[label=\alph*)]
         \item $|\PP(X_n(x) \in \cdot) - \PP(X_n(0) \in \cdot) | \le \PP(T_n \le x)$; \label{p1:a}
         \item Under \ref{a1}, $X_n(0) \eqd \widetilde X_n$ for $n \ge 0$; \label{p1:b}
         \item 
         Further, $\widetilde X_n \nearrow \widetilde X_\infty:= \max (0, \sup_{j\ge 1} [\tilde s_j - \sum_{k=0}^{j-1} \tilde t_k])$ a.s.~as $n \to \infty$. \label{p1:c}
      \end{enumerate}
   \end{proposition}

   \begin{proof}
      We note that the virtual workload for $G/G/\infty$ is the maximum of the workload associated with customer 0 (having remaining work $\max(x - T_j, 0)$ at time $T_j$) and the workload of all the customers to subsequently arrive to the system (given by $X_n(0)$ at time $T_n$). So, on $\{T_n > x\}$, $X_n(x) = X_n(0)$. Hence,
      \begin{align*}
         &\PP(X_n(x) \in A) - \PP(X_n(0) \in A) \\
         \le~&\PP(X_n(x) \in A, T_n > x) - \PP(X_n(0) \in A, T_n > x) \\
         &\qquad + \PP(T_n \le x) \\
         =~&\PP(T_n \le x),
      \end{align*}
      and similarly,
      \[ \PP(X_n(0) \in A) - \PP(X_n(x) \in A) \le \PP(T_n \le x), \]
      proving \ref{p1:a}.

      As for \ref{p1:b},
      \begin{align*}
         X_n(0) &= \max_{1 \le j \le n} [s_j - \sum_{k=j+1}^n t_k] \\
         &\eqd \max_{-n \le j \le -1} [s_j - \sum_{k=j+1}^{-1} t_k] \\
         &= \max_{1 \le j \le n} [\tilde s_j - \sum_{\ell=1}^{j-1} \tilde t_\ell] = \widetilde X_n,
      \end{align*}
      proving \ref{p1:b}, where the second equality follows from \ref{a1}.
      Finally, $(\widetilde X_n : n \ge 0)$ is a non-decreasing sequence, and hence $\widetilde X_n \nearrow \widetilde X_{\infty}$ a.s.~as $n \to \infty$.
   \end{proof}

   Note that Proposition \ref{p1} implies that $X_n(x) \Rightarrow \widetilde X_\infty$ as $n \to \infty$, where $\Rightarrow$ denotes weak convergence (on $[0, \infty]$). Set
   \[ X_j^* = \sup_{k \ge 0} [s_{j-k} - \sum_{\ell = j-k+1}^j t_\ell].\]
   Clearly, $X_j^* \eqd \widetilde X_\infty$ for $j \in \mathbb Z$.

   \begin{theorem}
      \label{thm1}
      Assume \ref{a1} and suppose that $\widetilde X_\infty$ is a finite-valued rv. Then,
      \begin{enumerate}[label=\alph*)]
         \item $((s_j, t_j, X_j^*) : -\infty < j < \infty)$ is a stationary sequence, and $(X_j^*:-\infty<j<\infty)$ satisfies the recursion \[ X_{j+1}^* = \max(X_j^* - t_{j+1}, s_{j+1})\]
         (so that it is a stationary solution of \eqref{SRS1}). \label{thm1:a}
         \item If $(X_j' : -\infty<j<\infty)$ is a stationary solution of \eqref{SRS1}, then $X_j' = X_j^*$ for $j \in \mathbb Z$. \label{thm1:b}
      \end{enumerate}
   \end{theorem}

   \begin{proof}
      For \ref{thm1:a}, it is clear that
      \begin{align*}
         X_{j+1}^* &= \sup_{k \ge 0} \left(s_{j+1-k} - \sum_{\ell = j - k + 2}^{j+1} t_\ell\right) \\
         &= \max \left(s_{j+1},\ \sup_{k \ge 1} [ s_{j+1-k} - \sum_{\ell=j-k+2}^j t_\ell - t_{j+1}]\right) \\
         &= \max\left(s_{j+1},\ \sup_{n \ge 0} [s_{j-n} - \sum_{\ell = j - n + 1}^j t_\ell] - t_{j+1}\right) \\
         &= \max\left(s_{j+1},\ X_j^* - t_{j+1}\right).
      \end{align*}
      For \ref{thm1:b}, we observe that for $n \ge 1$, $X_{-n}' \ge s_{-n}$, so
      \begin{align*}
         X_0' &= \max_{-n < j \le 0} \left(X_{-n}' - \sum_{\ell=-n+1}^0 t_\ell,\ s_j - \sum_{\ell=j+1}^0 t_\ell \right) \\
         &\ge \max_{-n < j \le 0} \left( s_{-n} -\sum_{\ell=-n + 1}^0 t_\ell,\ s_j - \sum_{\ell = j+1}^0 t_\ell \right) \\
         &= \max_{-n \le j \le 0} \left(s_j - \sum_{\ell = j+1}^0 t_\ell \right).
      \end{align*}
      Sending $n \to \infty$, we find that $X_0' \ge X_0^*$. We note that when $\sum_{\ell = -n +1 }^0 t_j > X_{-n}',$ then
      \[ X_0' = \max_{-n < j \le 0} \left( s_j - \sum_{\ell = j + 1}^0 t_\ell \right) \le X_0^*.\]
      So, if $f : \mathbb R_+ \to \mathbb R_+$ is a non-decreasing bounded function,
      \begin{align}
         \E f(X_0') &\le \E f(X_0')\,I\left( \sum_{\ell=0}^{n-1} \tilde t_\ell > X_{-n}' \right) \nonumber \\
         &\qquad \quad + \|f\|_\infty\,\PP \left(\sum_{\ell=0}^{n-1} \tilde t_\ell \le X_{-n}'\right) \nonumber \\
         &\le \E f(X_0^*) + \|f\|_\infty\,\PP \left(\sum_{\ell=1}^n t_\ell \le X_0'\right), \label{e21}
      \end{align}
      where we use the fact that $(X_{-n}' : -\infty < n < \infty)$ is a stationary solution of \eqref{SRS1} for the second inequality (and $\|f\|_\infty \triangleq \sup\{|f(x)| : x \in \mathbb R_+\}$). Since $X_0'$ is finite-valued, $P(\sum_{\ell=1}^n t_\ell \le X_0') \to 0$ as $n \to \infty$, and consequently \eqref{e21} implies that $\E f(X_0') \le \E f(X_0^*)$ for all bounded non-decreasing $f$.
      Since $X_0' \ge X_0^*$, it follows that $X_0' = X_0^*$. Similarly, $X_n' = X_n^*$ for $n \in \mathbb Z$.
   \end{proof}

   \begin{remark}
      Theorem \ref{thm1} is, of course, the $G/G/\infty$ analog of the Loynes construction for the $G/G/1$ queue. The analog to $\widetilde X_\infty$ is the rv $\max_{j \ge 0} \sum_{k=1}^j (\tilde s_{k-1} - \tilde t_k)$; see R.M.~Loynes (1962).
   \end{remark}

   We next turn to the question of when $\widetilde X_\infty$ is finite-valued. For this purpose, we strengthen \ref{a1} to:
   \begin{assumption}
      \label{a2}
      $(s_j : -\infty < j < \infty)$ is an independent and identically distributed (iid) sequence of positive rv's, independent of the stationary ergodic sequence $(t_j : -\infty < j < \infty)$ of positive rv's.
   \end{assumption}
   \noindent Because of the independence within the service time sequence, we refer to this model as the $G/GI/\infty$ queue.
   We put $\widetilde T_n = \sum_{k=1}^n \tilde t_k$ for $n \ge 0$ and let $F(\cdot) = P(s_1 \le \cdot)$ and $\bar F(\cdot) = P(s_1 > \cdot)$.

   \begin{proposition}
      \label{p2}~
      \begin{enumerate}[label=\alph*)]
         \item Assume \ref{a1} and $\E\,t_1 < \infty$. Then, $\widetilde X_\infty$ is finite-valued a.s.~if $\E\,s_1 < \infty$. \label{p2:a}
         \item Assume \ref{a2} and $\E\,t_1 < \infty$. Then, $\E\,s_1 < \infty$ if $\widetilde X_\infty$ is finite-valued a.s. \label{p2:b}
         \item Assume \ref{a2}. Then, $\widetilde X_\infty$ is finite-valued a.s.~if and only if $\PP (\sum_{k=1}^\infty \bar F(\widetilde T_k) < \infty) = 1$. \label{p2:c}
      \end{enumerate}
   \end{proposition}

   \begin{proof}
      If $\E\,s_1 < \infty$, $\sum_{k=1}^\infty \bar F(k) = \sum_{k=1}^\infty \PP (\tilde s_k > k) < \infty$, so the Borel-Cantelli lemma implies that $\tilde s_k = o(k)$ a.s., where $o(a_k)$ is a sequence with the property that $o(a_k) / a_k \to 0$ as $k \to \infty$. On the other hand, the ergodic theorem and stationarity imply that $\widetilde T_n / n \to \E [\tilde t_1 \mid \mathcal J] > 0$ a.s., where $\mathcal J$ is the invariant $\sigma$-algebra associated with $(\tilde t_j : -\infty < j < \infty)$. (The positivity of $\E[\tilde t_1 \mid \mathcal J]$ follows from the positivity of $\tilde t_1$.) Consequently, $\tilde s_n - \widetilde T_{n-1} \to -\infty$ a.s.~as $n \to \infty$, so that
      \[ \widetilde X_\infty = \max_{n \ge 1}\ [\tilde s_n - \widetilde T_{n-1}] < \infty \quad \text{a.s.}, \]
      proving \ref{p2:a}.

      If $\E\,s_1 = \infty$, then for each $m \in \mathbb Z_+$, $\sum_{k=1}^\infty \PP(\tilde s_k > m k) = \infty$, so $\olim_{k \to \infty} \tilde s_k / k \ge m$ a.s.~by the converse to the Borel-Cantelli lemma (due to the independence of the $\tilde s_k$'s under \ref{a2}). So, $\olim_{k \to \infty} \tilde s_k / k = \infty$ a.s., and it follows that $\widetilde X_\infty = \infty$ a.s., yielding \ref{p2:b}.

      For \ref{p2:c}, note that if $\widetilde X_\infty < \infty$ a.s., then there exists $\tilde m < \infty$ such that $\PP(\widetilde X_\infty > \tilde m) < 1$. So, $\sum_{k=1}^\infty I(\tilde s_k - \widetilde T_{k-1} > \tilde m) < \infty$ with positive probability. The event $\{ \sum_{k=1}^\infty I(\tilde s_k - \widetilde T_{k-1} > \tilde m) < \infty \} = \{ \sum_{k=1}^\infty \PP (\tilde s_k > \widetilde T_{k-1} + \tilde m \mid \mathcal G) < \infty\}$ by the Borel-Cantelli converse for independent rv's, where $\mathcal G = \sigma(\widetilde T_j : - \infty < j < \infty)$. Hence, $\PP(\sum_{k=1}^\infty \bar F(\widetilde T_{k-1} + \tilde m) < \infty) > 0$. The sequence $(\tilde t_n : -\infty < n < \infty)$ is an ergodic stationary sequence, because ergodicity and stationarity are preserved under time-reversal; see, for example, p.~458 of J.L.~Doob (1953).
      Furthermore, on $\{\widetilde T_i > \tilde m\}$, $\sum_{k=1}^\infty \bar F(\widetilde T_{i+k}) \le \sum_{k=1}^\infty \bar F(\widetilde T_{i+k} - \widetilde T_i + \tilde m)$, and hence $\PP(\sum_{i=1}^\infty \bar F(\widetilde T_i) < \infty) \ge \PP(\sum_{k=1}^\infty \bar F(\widetilde T_{k-1} + \tilde m) < \infty) > 0$. If we write $\sum_{i=1}^\infty \bar F (\widetilde T_i) = g(\tilde t_1, \tilde t_2, \dots)$, then $g(\tilde t_2, \tilde t_3, \dots) = \sum_{i=1}^\infty \bar F(\widetilde T_{i+1} - \widetilde T_1)$. Note that
      \begin{align*}
         &\E \sum_{i=1}^\infty \bar F(\widetilde T_i - \widetilde T_1) - \bar F(\widetilde T_i) \\
         =&~\sum_{i=1}^\infty \PP(\sum_{j=2}^i \tilde t_j < s_1 \le \sum_{j=1}^i \tilde t_j) \\
         =&~\sum_{i=1}^\infty \PP ( \sum_{j=1}^{i-1} t_j < s_1 \le \sum_{j=1}^i t_j) \\
         =&~\E \sum_{i=1}^\infty I(T_{i-1} < s_1 \le T_i)
      \end{align*}

      But $\sum_{i=1}^\infty I(T_{i-1} < s < T_i) = 1$, so $\E \sum_{i=1}^\infty (\bar F(\widetilde T_i - \widetilde T_{i-1}) - \bar F(T_i)) < \infty$. It follows that
      \[ I(\sum_{i=1}^\infty \bar F(\widetilde T_i) < \infty) = I(\sum_{i=1}^\infty \bar F(\widetilde T_{i+1} - \widetilde T_1) < \infty) \quad \text{a.s.} \]
      so that the event $\{ \sum_{i=1}^\infty \bar F(\widetilde T_i) < \infty\}$ is an invariant event. Consequently, the ergodicity of $(\tilde t_n : - \infty < n < \infty)$ implies that $\PP(\sum_{i=1}^\infty \bar F(\widetilde T_i) < \infty)$ is either 0 or 1. Since we have already shown this probability is positive, we may conclude that $\PP(\sum_{i=1}^\infty \bar F(\widetilde T_i) < \infty) = 1$.

      Conversely, if $\PP(\sum_{i=1}^\infty \bar F(\widetilde T_i) < \infty) = 1$, then $\sum_{i=1}^{\infty} \PP(\tilde s_{i+1} > \widetilde T_i \mid \mathcal G) < \infty$ a.s., so that the Borel-Cantelli lemma, applied conditional on $\mathcal G$, implies that $\sum_{i=1}^\infty I(\tilde s_{i+1} < \widetilde T_i) < \infty$ a.s., so that $\widetilde X_\infty < \infty$ a.s.
   \end{proof}

\section{Markov Recurrence Structure of the $GI/GI/\infty$ Queue}

In this section, we wish to study $X = (X_n : n \ge 0)$ as a Markov chain on the state space $S = [0, \infty]$. To this end, we now enhance our assumption on the model to:
\begin{assumption} \label{a3}
   $(s_n : - \infty < n < \infty)$ is an iid sequence of positive rv's independent of the iid sequence $(t_n : -\infty < n < \infty)$ of positive rv's.
\end{assumption}
\noindent Because of the independence on both the inter-arrival and service times, we denote this model as the $GI / GI / \infty$ queue. With this assumption in place and $X_0$ chosen independently of $((s_j, t_j) : j \ge 1)$, $X$ is a Markov chain. When $\PP(s_1 \le \bar s) = 1$ for some $\bar s < \infty$, the interval $[0, \bar s]$ is an absorbing set for $X$. Put $\PP_x(\cdot) \triangleq \PP(\cdot \mid X_0 = x)$ for $x \ge 0$.

This Markov chain contains embedded regenerative structure. It follows that when $X$ is recurrent, it will be Harris recurrent.
One vehicle for constructing regenerations is the following. Suppose we choose $m_0 \ge 1$ and $w_0 < \infty$ so that $\PP(s_1 \le w_0) \ge 1/2$ and $\PP(T_{m_0} > w_0) \ge 1/2$. At a time $\nu$ at which $X_\nu \le w_0$, all the work that has arrived to the system by time $T_\nu$ will have drained by time $T_\nu + w_0$. So, on the event $\{T_{\nu + m_0} - T_\nu > w_0\}$, $X_{\nu + m_0}$ will depend only on $((s_{\nu + i}, t_{\nu + i}): 1 \le i \le m_0)$, and
\begin{equation*}
   \begin{split}
      & \PP_x(X_{\nu + m_0} \in \cdot \mid X_j : 0 \le j \le \nu,~T_{\nu + m_0} - T_\nu > w_0) \\
      &= \PP_0(X_{m_0} \in \cdot \mid T_{m_0} > w_0) \triangleq \varphi(\cdot),
   \end{split}
\end{equation*}
so that $X_{\nu + m_0}$ then has a distribution $\varphi$ independent of $X_\nu$.

To construct an infinite sequence of such regeneration times, we let
\[ \tau_1 = \inf\{ n m_0 : X_{(n-1)m_0} \le w_0,~ T_{nm_0} - T_{(n-1)m_0}> w_0\},\]
and for $i \ge 1$, set
\[ \tau_{i+1} = \inf\{nm_0 > \tau_i : X_{(n-1) m_0} \le w_0,~T_{n m_0} - T_{(n-1) m_0} > w_0 \}. \]
In other words, the $\tau_i$'s are the times at multiples of $m_0$ at which $X_{\tau_i - m_0} \le w_0$ and $T_{\tau_i} - T_{\tau_i - m_0} > w_0$, so that the $X_{\tau_i}$'s then have distribution $\varphi$.

For $n \ge m_0$, let $Y_n = I(X_{n-m_0} \le w_0,~T_n - T_{n-m_0} > w_0)$, and set $R_i = Y_{i m_0}$ for $i \ge 1$. Then, $R_i$ is 1 or 0 depending on whether or not a regeneration occurs at time $T_{im_0}$. 
Conditionally on the event $\{T_{m_0}>w_0\}$, $\tau_1$ has the distribution of a typical regenerative cycle length. Put $\PP_\varphi(\cdot) = \int_S \varphi(dx) \PP_x(\cdot)$. Then, the number of regenerations is finite a.s.~when $\PP_\varphi(\tau_1 = \infty) > 0$. Set $u_i = \E_\varphi R_i$ with $u_0 = 1$, so that $(u_i : i \ge 0)$ is the \textit{renewal sequence} associated with the $\tau_i$'s.
(Conditional on $\tau_i < \infty$, $\tau_1$, $\tau_2 - \tau_1$, \dots, $\tau_i - \tau_{i-1}$ are iid positive rv's under $\PP_\varphi$.)
A well known consequence of renewal theory is the following result; see XIII.3 of W.~Feller (1968) for the proof.

\begin{proposition}
   \label{prop3}~
   \begin{enumerate}[label=\alph*)]
      
      \item $\PP_\varphi(\tau_1 < \infty) = 1$ if and only if $\sum_{n=0}^\infty u_n = \infty$. \label{prop3:a}
      \item $\E_\varphi\,\tau_1 < \infty$ if and only if $\olim_{n \to \infty} (1/n) \sum_{i=0}^{n-1} u_i > 0$.
   \end{enumerate}
\end{proposition}

We will utilize Proposition \ref{prop3} to prove our main theorem on recurrence classification for the Markov chain $X$.

\begin{theorem} \label{thm2}
   Assume \ref{a3}.
   \begin{enumerate}[label=\alph*)]
      \item If $\sum_{n=1}^\infty \E \prod_{i=1}^n F(y + \widetilde T_i) < \infty$, then for $x \in \mathbb R_+$, \[ \E_x \sum_{n=0}^\infty I(X_n \le y) < \infty.\] \label{thm2:a}
      \item If $\sum_{n=1}^\infty \E \prod_{i=0}^{n-1} F(w_0 + \widetilde T_i) = \infty$, then $X$ is a Harris recurrent Markov chain on $\mathbb R_+$. \label{thm2:b}
      \item $\PP(\sum_{i=1}^\infty \bar F(\widetilde T_i) < \infty) = 1$ if and only if $X$ is a positive recurrent Harris chain on $\mathbb R_+$ with stationary distribution $\pi$ given by $\pi(\cdot) = \PP(\widetilde X_\infty \in \cdot)$. \label{thm2:c}
   \end{enumerate}
\end{theorem}
\begin{proof}
   We note that
   \begin{align*}
      \PP_x(X_n \le y) &= \PP(x - \sum_{i=1}^n t_i \le y,~s_1 - \sum_{i=2}^n t_i \le y, \dots, s_{n-1} - t_n \le y,~s_n \le y) \displaybreak[1] \\
      &= \E\,I(\sum_{i=1}^n \tilde t_i \ge x - y,~\tilde s_{n-1} \le y + \sum_{i=1}^{n-1} \tilde t_i, \dots, \tilde s_2 \le y+\tilde t_1,~\tilde s_1 \le y) \displaybreak[1]\\
      &= \E\,I(\widetilde T_n \ge x-y) \prod_{i=1}^n F(y + \widetilde T_{i-1}) \displaybreak[1]\\
      &\le \E \prod_{i=1}^{n-1} F(y + \widetilde T_i) \displaybreak[1]\\
      &= \E \exp(\sum_{i=1}^{n-1} \log(1 - \bar F(y + \widetilde T_i))) \displaybreak[1]\\
      &\le \E \exp(-\sum_{i=1}^{n-1} \bar F(y + \widetilde T_i)).
   \end{align*}
So,
\[ \E_x \sum_{n=0}^\infty I(X_n \le y) \le \sum_{n=0}^\infty \E \exp(-\sum_{i=1}^{n-1} \bar F(y + \widetilde T_i)) < \infty;\]
proving \ref{thm2:a}.

For \ref{thm2:b}, we choose $x_0$ so that $\PP_\varphi(X_0 \le x_0) > 0$. Then,
\begin{align*}
   \PP_\varphi(X_n \le w_0) &\ge \PP_\varphi(X_n \le w_0,\,X_0 \le x_0) \displaybreak[1] \\
   &\ge \PP(X_0 - \sum_{i=1}^n t_i \le w_0,~s_1 - \sum_{i=2}^n t_i \le w_0, \dots, s_n \le w_0,~X_0 \le x_0) \displaybreak[1]\\
   \begin{split}
      &\ge \PP(\widetilde T_n \ge x_0 - w_0,~\tilde s_1 \le \widetilde T_{n-1} + w_0,~\tilde s_2 \le \widetilde T_{n-2} + w_0, \dots, \tilde s_{n-1} + \widetilde T_1 \le w_0) \\
      &\qquad \cdot\PP_\varphi(X_0 \le x_0)
   \end{split} \displaybreak[1]\\
   &\ge \E \prod_{i=0}^{n-1} F(\widetilde T_i + w_0) \PP_\varphi(X_0 \le x_0) - \PP(\widetilde T_n < x_0 - w_0). \displaybreak[1]
\end{align*}
It is well known that $\PP(\widetilde T_n < x_0 - w_0)$ decays to 0 geometrically in $n$ and thus $\sum_{n=0}^\infty \PP(\widetilde T_n < x_0 - w_0) < \infty$ (see, for example, W. Feller (1971), p.~185). Consequently, we conclude in the presence of our assumption in \ref{thm2:b} that
\begin{align*}
   \infty &= \sum_{n=1}^\infty \PP_\varphi (X_{nm_0} \le w_0) \PP(T_{m_0} > w_0) \\
   &= \sum_{n=2}^\infty \E_\varphi R_n = \sum_{n=2}^\infty u_n,
\end{align*}
proving that the $\tau_n$'s comprise a recurrent regenerative sequence. Clearly, $A \triangleq [0, w_0]$ is a small set for $X$ (since $\PP_x(X_{m_0} \in \cdot) \ge \PP(T_m > w_0) \varphi(\cdot)$ for $x \in A$), and the regenerative recurrence implies that $\PP_\varphi(X_{nm_0} \in A~\text{infinitely often}) = 1.$ Hence, if $B \triangleq \{X_{n m_0} \in A~\text{infinitely often}\}$ and $h(x) \triangleq \E_x I(B)$, $\E_\varphi h(X_{km_0}) = 1$ for $k \ge 1$. Then, part \ref{p1:a} of Proposition \ref{p1} implies that
\[ \E_x h(X_{k m_0}) - \E_\varphi h(X_{k m_0}) \to 0\]
as $k \to \infty$, so that $\E_x h(X_{km_0}) \to 1$ as $k \to \infty$. But $h(X_n)$ is a $\PP_x$-martingale adapted to $(\mathcal F_n : n \ge 0)$, where $\mathcal F_n = \sigma(X_j : 0 \le j \le n)$, so $h(x) = \E_x h(X_{k m_0})$, and hence $h(x) = 1$ for $x \in \mathbb R_+$, proving that $\PP_x(X_{n m_0} \in A~\text{infinitely often}) = 1$ for $x \in S$. Therefore, $X$ is a Harris recurrent Markov chain.

As for \ref{thm2:c}, we have already shown in part \ref{p2:c} of Proposition \ref{p2} that $X$ has a stationary distribution $\pi$ if and only if $\PP(\sum_{k=1}^\infty \bar F(\widetilde T_k) < \infty) = 1$. It remains only to verify that this condition implies that $X$ is Harris recurrent.

Note that $z \in [0, \bar F(w_0)]$, $\log(1-z) \ge -c_0 z$ for some positive $c_0$. Hence,
\begin{align*}
   \E \prod_{i=1}^n F(\widetilde T_i + w_0) &= \E \exp(\sum_{i=1}^n \log (1-\bar F(\widetilde T_i + w_0))) \displaybreak[1]\\
   &\ge \E \exp(-c_0 \sum_{i=1}^n \bar F(\widetilde T_i + w_0)) \displaybreak[1]\\
   &\ge \E \exp(-c_0 \sum_{i=1}^n \bar F(\widetilde T_i + \sum_{j=-m_0}^{-1} \tilde t_j)) I(\sum_{j=-m_0}^{-1} \tilde t_j > w_0) \displaybreak[1]\\
   &= \E \exp(-c_0 \sum_{i=1}^{m_0 + n} \bar F(\widetilde T_i)) I(\widetilde T_{m_0} > w_0).
\end{align*}
Since the Bounded Convergence Theorem implies that
\[ \E \exp(-c_0 \sum_{i=1}^{m_0 + n} \bar F(\widetilde T_i)) I(\widetilde T_{m_0} > w_0) \searrow \E \exp(-c_0 \sum_{i=1}^\infty \bar F(\widetilde T_i)) I(\widetilde T_{m_0} > w_0) > 0,\]
under the assumption that $\sum_{i=1}^\infty \bar F(\widetilde T_i) < \infty$ a.s., it follows that
\[ \sum_{n=1}^\infty \E \prod_{i=1}^n \bar F(\widetilde T_i + w_0) = \infty,\]
so that \ref{thm2:b} yields the Harris recurrence.
\end{proof}

\begin{remark}
   Note that under \ref{a3}, $(\widetilde T_i : i \ge 0) \eqd (T_i : i \ge 0)$, so the $T_i$'s may replace the $\widetilde T_i$'s in the statement of Theorem \ref{thm2}.
\end{remark}

Our next result provides simplified sufficient conditions for \ref{thm2:a} and \ref{thm2:b} of Theorem \ref{thm2}.

\begin{proposition} \label{p4}
   Assume \ref{a3}.
   \begin{enumerate}[label=\alph*)]
      \item If \[ \sum_{n=1}^\infty \E \exp(-\sum_{i=1}^n \bar F(T_i + w_0)) < \infty, \]
      then $X$ is transient. \label{p4:a}
      \item If there exists $c > 1$ such that \[ \sum_{n=1}^\infty \E \exp(-c \sum_{i=1}^n \bar F(T_i + w_0)) = \infty,\]
      then $X$ is a Harris recurrent Markov chain. \label{p4:b}
   \end{enumerate}
\end{proposition}
\begin{proof}
   Note that since $\log(1-x) \le -x$ for $x \ge 0$,
   \[ \E \prod_{i=1}^n F(T_i + w_0) \le \E \exp(-\sum_{i=1}^n \bar F(T_i + w_0)),\]
   proving \ref{thm2:a}.

   For \ref{thm2:b}, let $N(t) = \max\{ n \ge 0 : T_n \le t\}$ for $t \ge 0$ be the renewal counting process corresponding to the $T_i$'s. Select $z_0$ so large that $\log(1 - \bar F(z)) \ge -c \bar F(z)$ for $z \ge z_0$. Then,
   \begin{align*}
      &\E \sum_{n=1}^\infty \prod_{i=0}^{n-1} F(T_i + w_0) \displaybreak[1]\\
      &\ge \E \sum_{n=N(z_0) + 2}^\infty \prod_{i=0}^{N(z_0) + 1} F(T_i + w_0) \exp(-c \sum_{i=N(z_0) + 2}^n \bar F(T_i + w_0)) \displaybreak[1]\\
      &\ge \E \sum_{n=1}^\infty \prod_{i=0}^{N(z_0) + 1} F(T_i + w_0) \exp(-c \sum_{i=N(z_0) + 2}^{N(z_0) + 1 + n} \bar F (\sum_{j=N(z_0) + 2}^i t_j + w_0)) \displaybreak[1]\\
      &= \E \prod_{i=0}^{N(z_0) + 1} F(T_i + w_0) \sum_{n=1}^\infty \E [\exp(-c \sum_{i=N(z_0) + 2}^{N(z_0) + 1 + n} \bar F(\sum_{j=N(z_0)+2}^i t_j + w_0)) \mid \mathcal F_{N(z_0) + 1}] \displaybreak[1]\\
      &\ge \E F(w_0)^{N(z_0) + 1} \sum_{n=1}^\infty \E \exp(-c \sum_{i=1}^n \bar F(T_i + w_0)), \displaybreak[1]
   \end{align*}
   so divergence of the latter sum implies Harris recurrence of $X$.
\end{proof}

We can now apply our recurrence theory to specific examples.

\begin{example}
   Suppose that
   \begin{equation} \label{e31}
      P(s_1 > x) \sim d_1 x^{-\alpha}
   \end{equation}
   and
   \begin{equation} \label{e32}
      P(t_1 > x) \sim d_2 x^{-\beta}
   \end{equation}
   as $x \to \infty$, for $d_1, d_2, \alpha, \beta > 0$ and $\alpha, \beta < 2$, where we write $a(x) \sim b(x)$ as $x \to \infty$ whenever $a(x) / b(x) \to 1$ as $x \to \infty$. If \eqref{e32} is strengthened to
   \[ P(t_1 > x) = d_2 x^{-\beta} + d_3 x^{-\delta} + O(x^{-\gamma}) \]
   as $x \to \infty$ for $\beta < \delta < \gamma$, where $O(x^{-\gamma})$ is monotone decreasing (and $O(a(x))$ denotes a function for which $O(a(x))/a(x)$ remains bounded as $x \to \infty$), then J.~Chover (1966) established the law of the iterated logarithm
   \[ \olim_{n \to \infty} \left(\frac{T_n}{n^{1/\beta}} \right)^{\frac{1}{\log \log n}} = e^{-1/\beta} \quad \text{a.s.}\]
   as $n \to \infty$. It follows that for $\epsilon > 0$, $\sum_{n=0}^\infty \bar F(T_n)$ converges a.s.~when
   \[\sum_{n=0}^\infty ((e^{-1/\beta} - \epsilon)^{\log \log n} n^{1/\beta})^{-\alpha} < \infty,\]
   which arises when $\alpha > \beta$. So, the $GI/GI/\infty$ queue is a positive recurrent Harris chain in this setting. When $\alpha < \beta$, $\sum_{n=0}^\infty \bar F(T_n) = \infty$ a.s., and $X$ is not positive recurrent.
\end{example}

\begin{example}
   If $T_n = r n$ for $r > 0$ (so that the arrival process is deterministic) and $s_1$ satisfies \eqref{e31}, then
   \[ \E \exp(-\sum_{k=1}^n \bar F(T_k + w_0)) = \exp(-\sum_{k=1}^n \bar F(rn + w_0)).\]
   When $\alpha < 1$, $X$ is then transient by part \ref{p4:a} of Proposition \ref{p4}, whereas if $\alpha > 1$, then $\E\,s_1 < \infty$, so $X$ is positive recurrent. If $\PP(s_1 > x) = c x^{-1}$ for $x$ sufficiently large, then $\sum_{n=1}^\infty \exp(-\sum_{k=1}^n \bar F(rk+w_0))$ diverges or converges in accordance with
   \begin{align*}
      & \sum_{n=1}^\infty \exp(-d_1 \sum_{k=1}^n (rk + w_0)^{-1}) \\
      &= \sum_{n=1}^\infty \exp\left(- \frac{d_1}{r} \log \left(r + \frac{w_0}{r}\right) + O(1)\right) \\
      &= \exp(O(1)) \sum_{n=1}^\infty n^{-d_1/r}.
   \end{align*}
   Consequently, when $\alpha = 1$, $X$ is transient when $d_1 > r$ and null recurrent when $d_1 < r$, according to Proposition \ref{p4}. At $d_1 = r$, we apply part \ref{thm2:b} of Theorem \ref{thm2} to conclude that $X$ is null recurrent.
\end{example}

In this $GI/GI/\infty$ setting, there is a well-developed parallel theory for recurrence and transience of the corresponding $GI/GI/1$ queue (i.e.~under \ref{a3}). When at least one of $\E\,s_1$ and $\E\,t_1$ is finite, it is well known that $(W_n : n \ge 0)$ is a positive recurrent Harris chain that hits state 0 infinitely often if and only if $\E\,s_1 < \E\,t_1$. The chain is null recurrent when $\E\,s_1 = \E\,t_1$, and transient when $\E\,s_1 > \E\,t_1$.

The mathematically delicate case arises when $\E\,s_1 = \E\,t_1 = \infty$. In this setting, recurrence is intimately connected to subtle properties of the \textit{random walk} $\Gamma = (\Gamma_n : n \ge 0)$, where $\Gamma_n = \sum_{j=1}^n (s_j - t_{j+1})$; and to its maximum $M_\infty \triangleq \max_{n \ge 0} \Gamma_n$. The chain $(W_n : n \ge 0)$ is positive recurrent when $M_\infty < \infty$ a.s.~and null recurrent when $\olim_{n \to \infty} \Gamma_n = \infty$ a.s.~and $\ulim_{n \to \infty} \Gamma_n = -\infty$. General classification results in terms of $(\Gamma_n : n \ge 0)$ are classical (see, for example, F.~Spitzer (1964)), but they are difficult to apply. K.B.~Erickson (1973) made a significant contribution by providing a characterization of the recurrence / transience theory, expressed purely in terms of the marginal distribution of $\xi_1 = s_0 - t_1$. Put $\xi_1^+ = \max(\xi_1, 0)$, $\xi_1^- = \max(-\xi_1, 0)$. The following is due to Erickson.

\begin{theorem}
   Assume $\PP(\xi_1 > 0)$ and $\PP(\xi_1 < 0) > 0$. Let 
   \[ \hat J_+ = \int_0^\infty \frac{x}{\E \min(\xi_1^-, x)} \PP(\xi_1^+ \in dx), \quad \hat J_1^- =  \int_0^\infty \frac{x}{\E \min(\xi_1^+, x)} \PP(\xi_1^- \in dx).\]
   If $\E\,|\xi_1|=\infty$, then $\hat J_++\hat J_- = \infty$, and 
   \begin{enumerate}[label=\alph*)]
      \item $\lim_{n \to \infty} \Gamma_n / n = \infty$ a.s.~if and only if $\hat J_- < \infty$;
      \item $\lim_{n \to \infty} \Gamma_n / n = -\infty$ a.s.~if and only if $\hat J_+ < \infty$;
      \item $\olim_{n \to \infty} \Gamma_n / n = -\ulim_{n \to \infty} \Gamma_n / n = \infty$ a.s.~if and only if $\hat J_+ = \hat J_- = \infty$.
   \end{enumerate}
\end{theorem}

For our current purposes, it is convenient to reformulate these conditions. For $t \ge 0$, put $m_+(t) = \E\,\min(s_1,t)$ and $m_-(t)=\E\,\min(t_1,t)$. Then, $\hat J_+ < \infty$ if and only if 
\[ J_+ \triangleq \E\, \left(\frac{s_1}{m_-(s_1)}\right) < \infty,\]
whereas $\hat J_- < \infty$ if and only if 
\[ J_- \triangleq \E\,\left(\frac{t_1}{m_+(t_1)}\right) < \infty.\]

Recall the notation $N(t) = \max\{n\ge 0: T_n\leq t\}$. By Wald's identity (see, for example, p.~397 in W.~Feller (1971)), it is easily seen that 
\[ \frac{x}{m_-(x)} \le \E\,N(x) + 1 \le \frac{2x}{m_-(x)}, \]
so that $J_+ < \infty$ if and only if 
\begin{equation} \label{e3q}
   \E\,N(s_1) < \infty.
\end{equation}
Thus, if $\E |\xi|=\infty$, the $GI / GI / 1$ queue is positive recurrent if and only if \eqref{e3q} is in force. Note that our necessary and sufficient condition for positive recurrence of $GI/GI/\infty$ is closely related. In particular, $\sum_{n=0}^\infty \bar F(T_n) < \infty$ a.s.~is equivalent to 
\begin{equation} \label{e3qq}
   \E[N(s_1) \mid \mathcal G] < \infty \quad \text{a.s.} 
\end{equation}
So, $GI/GI/\infty$ requires the conditional expectation to be finite a.s., whereas $GI/GI/1$ requires the expectation to be finite.

Let us revisit Example 1. Assume in addition that $\alpha <1$ and $\beta < 1$. Then $J_+<\infty$ if and only if $\alpha > \beta$, so both $GI/GI/1$ and $GI/GI/\infty$ queues are positive Harris recurrent in this case.

\begin{remark}
 The discussion in this section provides a few settings in which sufficient conditions for recurrence/transience, expressed in terms of marginal distributions, can be developed. An open research question is the development of a more comprehensive theory, again expressed in terms of marginal distributions, for settling such recurrence/transience issues for the infinite server queue. 
C.M.~Goldie and R.A.~Maller (2000) provide one such framework in their Theorem 2.1 for a related stochastic recursion. 
\end{remark}

\section{Tail Asymptotics for the Distribution of $\widetilde X_\infty$}

We now wish to study the tail probability $\PP(\widetilde X_\infty > x)$ as $x \to \infty$.

\begin{theorem} \label{thm3}
   Assume \ref{a2}.
   \begin{enumerate}[label=\alph*)]
      \item If $\bar F(x) \sim \delta \exp(-\mu x)$ as $x \to \infty$ for $\delta, \mu > 0$, then
      \[ \PP(\widetilde X_\infty > x) \sim \delta \exp(-\mu x) \sum_{n=0}^\infty \E \exp(-\mu T_n) \]
      as $x \to \infty$. \label{thm3:a}
      \item If $\widetilde F(x) \sim \delta x^{-\alpha}$ as $x \to \infty$ for $\delta > 0$, $\alpha > 1$, and $\E\,t_1 < \infty$, then \[ \PP(\widetilde X_\infty > x) \sim 
      \frac{\delta}{\E\,t_1 (\alpha-1)} x^{1-\alpha}\]
      as $x \to \infty$. \label{thm3:b}
   \end{enumerate}
\end{theorem}

\begin{proof}
   We start by observing that
   \begin{align*}
      \PP(\widetilde X_\infty > x) &= 1 - \PP(\widetilde X_\infty \le x) \\
      &= 1 - \E \exp(\sum_{i=0}^\infty \log (1 - \bar F(x+T_i))) \\
      &= 1 - \E \exp(- (1+o(1))\sum_{i=0}^\infty \bar F(x + T_i) )
   \end{align*}
   where $o(1)$ is deterministic (and depends only on $x$). Since $y - y^2/2 \le 1 - \exp(-y) \le y$ for $y \ge 0$,
   \begin{align*}
      & 1 - \E \exp(-(1+o(1))\sum_{i=0}^\infty \bar F(x + T_i)) \displaybreak[1] \\
      &= (1+o(1))\E \sum_{i=0}^\infty \bar F(x + T_i) \displaybreak[1]\\
      &= (1+o(1))\E \sum_{i=0}^\infty \delta \exp(-\mu (x +  T_i)) \displaybreak[1]\\
      &= (1+o(1))\delta \exp(-\mu x) \sum_{i=0}^\infty \E \exp(-\mu  T_i), \displaybreak[1]
   \end{align*}
   proving \ref{thm3:a}.

   For \ref{thm3:b},
   \begin{align*}
      &1 - \E \exp(-\sum_{i=0}^\infty \bar F(x + 
      T_i)(1 + o(1))) \displaybreak[1] \\
      & = (1+o(1))\E \sum_{i=0}^\infty \delta (x + 
      T_i)^{-\alpha}  \displaybreak[1] \\
      &= (1+o(1))\delta x^{1-\alpha} \E \sum_{i=0}^\infty \frac{1}{x} \left(1 + \frac{
      T_i}{x}\right)^{-\alpha} 
      \displaybreak[1] \\
      &= (1+o(1))\delta x^{1-\alpha} \E \sum_{i=0}^\infty \frac{1}{x} \left( 1 + \frac{i}{x} \E\,t_1 (1 + o(1))\right)^{-\alpha} 
      \displaybreak[1] \\
      &= (1+o(1))\delta x^{1-\alpha} \int_0^\infty (1 + y\,\E\,t_1)^{-\alpha}\,dy~
      \displaybreak[1] \\
      &= \delta x^{1-\alpha} \frac{1}{\E\,t_1 (\alpha-1)} (1 + o(1)),
   \end{align*}
   providing the proof of \ref{thm3:b}.
\end{proof}

For $GI/GI/\infty$, the tail of $\widetilde X_\infty$ is, up to a constant, that of the service time distribution when the tail is asymptotically exponential, whereas the steady-state $M_\infty$ for $GI/GI/1$ exhibits Cram{\`e}r-Lundberg asymptotics in this setting; see S.~Asmussen (2003), XIII.5. On the other hand, for Pareto-type tails, the asymptotics for $GI/GI/\infty$ and $GI/GI/1$ are the same, up to a constant factor.\\

{\bf Acknowledgement:} The authors thank the reviewer for important comments and suggestions.



\begin{thebibliography}{99}\small

\bibitem{Alt2005} 
E.~Altman. 
 ``On stochastic recursive equations and infinite server queues.'' 
 In {\it Proceedings IEEE 24th Annual Joint Conference of the IEEE Computer and Communications Societies}, Vol. 2 (2005), 1295--1302.

 \bibitem{Asm2003} 
 S.~Asmussen. 
``Applied Probability and Queues.'' 
 Springer, 2003.

\bibitem{Bac1995} 
F.~Baccelli, S.~Foss. 
``On the saturation rule for the stability of queues.'' 
 {\it Journal of Applied Probability}, 32, 2 (1995), 494--507.

 \bibitem{Cho1966} 
 J.~Chover. 
 ``A law of the iterated logarithm for stable summands.'' 
 {\it Proceedings of the American Mathematical Society}, 17 (1966), 441--443.

 \bibitem{Com2002}  
 F.~Comets, R.~Fernandez, P.A.~Ferrari. 
 ``Processes with long memory: regenerative construction and perfect simulation.'' 
 {\it The Annals of Applied Probability}, 12, 2 (2002), 921--943.

 \bibitem{Doo1953} 
 J.L.~Doob. 
 ``Stochastic Processes.'' 
 John Wiley \& Sons, Inc., New York; Chapman \& Hall, Ltd., London, 1953.

 \bibitem{Eri1973}
 K.B.~Erickson. 
 ``The strong law of large numbers when the mean is undefined.'' 
 {\it Transactions of the American Mathematical Society}, 185 (1973), 371--381.
 
  \bibitem{Fel1968} 
  W.~Feller.  
 ``An Introduction to Probability Theory and Its Applications, Volume I.'' 
 John Wiley \& Sons, Inc., New York, 1968.

 \bibitem{Fel1971} 
 W.~Feller. 
 ``An Introduction to Probability Theory and Its Applications, Volume II.'' 
 John Wiley \& Sons, Inc., New York, 1971.
 
 \bibitem{Gol1991} 
 C.M.~Goldie. 
 ``Implicit renewal theory and tails of solutions of random equations.'' 
 {\it The Annals of Applied Probability}, 1, 1 (1991), 128--166. 
 
 \bibitem{Gol2000} 
 C.M.~Goldie, R.A.~Maller. 
 ``Stability of perpetuities.'' 
 {\it  The Annals of Probability}, 28, 3 (2000), 1195--1218. 


\bibitem{Loy1962}
R.M.~Loynes. 
``The stability of a queue with non-independent interarrival and service times.'' 
{\it Proceedings of the Cambridge Philosophical Society}, 58 (1962), 497--520.


\bibitem{Pop2025}
S.~Popov. 
``On transience of {${\rm M}/{\rm G}/\infty$} queues.'' 
 {\it Journal of Applied Probability}, 62, 2 (2025), 572--575.



 \bibitem{Spi1964} 
 F.~Spitzer. 
 ``Principles of Random Walk.''  
 Springer, 1964.

 \end{thebibliography}
\end{document}